\documentclass[amsthm]{elsart}
\usepackage{ifpdf}
\usepackage{graphicx,amssymb,lineno}
\usepackage{subfig}
\usepackage{longtable}
\usepackage{epstopdf}
\usepackage{float}
\usepackage{cite}
\usepackage{mathrsfs}
\usepackage{latexsym,lineno}
\usepackage{epsfig}
\usepackage{color}
\usepackage{amsmath}\usepackage{fleqn}\usepackage{verbatim}\usepackage{epsf}
\usepackage{amsthm}\usepackage{graphicx, float}\usepackage{graphicx}
\usepackage{amsfonts}\usepackage{amssymb}\usepackage{graphpap}
\usepackage{epic}\usepackage{curves}
\ifpdf
\usepackage[%
  pdftitle={Instructions for use of the document class
    elsart},%
  pdfauthor={Simon Pepping},%
  pdfsubject={The preprint document class elsart},%
  pdfkeywords={instructions for use, elsart, document class},%
  pdfstartview=FitH,%
  bookmarks=true,%
  bookmarksopen=true,%
  breaklinks=true,%
  colorlinks=true,%
  linkcolor=blue,anchorcolor=blue,%
  citecolor=blue,filecolor=blue,%
  menucolor=blue,pagecolor=blue,%
  urlcolor=blue]{hyperref}
\else
\usepackage[%
  breaklinks=true,%
  colorlinks=true,%
  linkcolor=blue,anchorcolor=blue,%
  citecolor=blue,filecolor=blue,%
  menucolor=blue,pagecolor=blue,%
  urlcolor=blue]{hyperref}
\fi

\makeatletter
\def\elsartstyle{%
    \def\normalsize{\@setfontsize\normalsize\@xiipt{14.5}}
    \def\small{\@setfontsize\small\@xipt{13.6}}
    \let\footnotesize=\small
    \def\large{\@setfontsize\large\@xivpt{18}}
    \def\Large{\@setfontsize\Large\@xviipt{22}}
    \skip\@mpfootins = 18\p@ \@plus 2\p@
    \normalsize
}
\@ifundefined{square}{}{}
\makeatother

\begin{document}

\begin{frontmatter}
\title{On the largest $A_{\alpha}$-spectral radius of cacti}

\journal{~~~~}

\author[SW]{Shaohui Wang}\ead{shaohuiwang@yahoo.com; wangs@savannahstate.edu},
\author[CW]{Chunxiang Wang}\ead{wcxiang@mail.ccnu.edu.cn},
\author[LIU]{Jia-Bao Liu}\ead{liujiabaoad@163.com},
\author[BW]{Bing Wei}\ead{bwei@olemiss.edu}

\address[SW]{
  Department of Mathematics, Savannah State University, Savannah, GA 31404, USA}
\address[CW]{ School of Mathematics and Statistics, Central China Normal University, Wuhan,
430079, P.R. China}
\address[LIU]{
  School of Mathematics and Physics, Anhui Jianzhu University, Hefei, 230601, P.R. China}
\address[BW]{
  Department of Mathematics, University of Mississippi, University, MS 38677, USA}

\corauth[cor]{Corresponding author: Chunxiang Wang.}

\begin{abstract}
Let $A(G)$ be the adjacent matrix and $D(G)$ the diagonal matrix of the degrees of a graph $G$, respectively. For $0 \leq \alpha \leq 1$,  
the $A_{\alpha}$ matrix  $A_{\alpha}(G) = \alpha D(G) +(1-\alpha)A(G)$ is given by Nikiforov. Clearly, $A_{0} (G)$ is the adjacent matrix and $2 A_{\frac{1}{2}}$ is the signless Laplacian matrix.  A cactus  is a connected graph such that  any two of its cycles have at most one common vertex, that is an extension of the tree. The $A_{\alpha}$-spectral radius of a cactus graph with $n$ vertices and $k$ cycles is explored. The outcomes obtained in this paper can imply previous bounds of Nikiforov et al., and Lov\'{a}sz and Pelik\'{a}n.  In addition,  the corresponding extremal graphs are determined. Furthermore, we proposed all eigenvalues of such extremal cacti. Our results extended and enriched previous known results.

\vskip 2mm \noindent {\bf Keywords:}  Signless Laplacian, adjacency matrix, tree, cacti. \\
\noindent{\bf AMS subject classification:} 05C50, 15A48.
 \end{abstract}

\end{frontmatter}

\section{Introduction}
We consider simple finite graph $G$ with vertex set $V(G)$ and edge set $E(G)$ throughout this work. The order of a graph is $ |V(G)| = n$ and the size is $|E(G)| = m$. For a vertex $v \in V(G)$, the neighborhood of  $v $ is the set $N(v) = N_G(v) = \{w \in V(G), vw \in E(G)\}$, and $d_G(v)$ (or briefly $d_v$)  denotes the degree of $v$ with $ d_{G}(v) = |N(v)|$. For $L \subseteq V(G)$ and  $R \subseteq E(G)$,   let $G[L]$ be the subgraph of $G$ induced by $L$, $G - L$ the subgraph induced by $V(G) - L$ and $G - R$  the subgraph of G obtained by deleting $R$.  Let $w(G - L)$ be the number of components of $G - L$, and  $L$ be a cut set if $w(G - L) \geq 2$.   If $e$ is an edge of $G$ and $w(G-e) \geq 2$, then $e$ is a cut edge of $G$. If $G-e$ contains at least two components, each of which contains at least two vertices, then $e$ is called a proper cut edge of $G$. Let $K_n, P_n$ and $S_n$ denote the clique, the path and the star on $n$ vertices, respectively.  If $P_k  = v_1v_2 \cdots v_k$ is a subgraph of $G$ and $v_1$ is a cut vertex of degree at least $3$, then $P_k$ is called a pendant path in $G$.

 Let $A(G)$ be the adjacency matrix and $D(G)$ the diagonal matrix of the degrees of $G$. The signless Laplacian matrix of $G$ is considered as $$Q(G) = D(G)+ A(G).$$ As the successful considerations on $A(G)$ and $Q(G)$,  Nikiforov \cite{005} proposed the matrix  $A_{\alpha}(G)$  of a graph $G$ $$A_{\alpha}(G) = \alpha D(G) +(1-\alpha)A(G),$$ for $\alpha \in [0,1]$.  It is not hard to see that if $\alpha =0, A_{\alpha}$ is the adjacent matrix, and if $\alpha = \frac{1}{2}$, then $2A_{\frac{1}{2}}$ is the signless Laplacian matrix of $G$. 

 In the mathematical literature there are numerous studies of properties of the (signless, $A_{\alpha}$)spectral radius. For instance, Chen \cite{004} explored properties of spectra of graphs and line graphs. 
 Lov\'{a}sz and  J. Pelik\'{a}n \cite{003} deduced the spectral radius of trees.   Cvetkovi\'{c} \cite{006} proposed the spectra of signless Laplacians of graphs and discuss  a related spectral theory of graphs. Zhou \cite{007} obtained the bounds of signless Laplacian spectral radius and its hamiltonicity. Lin and Zhou \cite{008} studied graphs with at most one signless Laplacian eigenvalue exceeding three.
 In addition to the thriving  considerations of the  spectral radius,  the  $A_{\alpha}$-spectral radius would be attractive. 
 
We first introduce some interesting properties for the  $A_{\alpha}$ matrix. 
Let $G$ be a graph with vertex set $V(G)= \{u_1,u_2, \cdots, u_n\}$ and edge set $E(G)$. The adjacent matrix of $G$ is $A(G)$, and the  $(i,j)$-entry of $ A(G)$ is $1$ if $u_iu_j \in E(G)$, and otherwise $0$. Denote  the eigenvalues of $A_{\alpha}(G)$ by $\lambda_1(A_{\alpha}(G)) \geq \lambda_2(A_{\alpha}(G)) \geq \cdots \geq \lambda_n(A_{\alpha}(G))$. The largest eigenvalue $\rho(G) := \lambda_1(A_{\alpha}(G))$ is defined as the $A_{\alpha}$-spectral radius of $G$. 
Denote by $X = (x_{u_1}, x_{u_2}, \cdots, x_{u_n})^T$ a real vector.  As $ A_{\alpha}(G) = \alpha D(G) + (1 - \alpha) A(G)$,  the quadratic form of $X^T A_{\alpha}(G)X$ can be written as  
\begin{eqnarray}
X^T A_{\alpha}(G)X = \alpha \sum_{u_i \in V(G)} x_{u_i}^2 d_{u_i} + 2(1 - \alpha) \sum_{u_iu_j \in E(G)} x_{u_i}x_{u_j}.
\label{01}
\end{eqnarray}
Because $A_{\alpha}(G)$ is a real symmetric matrix, and by Rayleigh principle, we have the important formula
\begin{eqnarray}
\rho(G) = max_{X \neq 0} \frac{X^T A_{\alpha}(G)X}{X^T X}.
\label{02}
\end{eqnarray} 
If $X$ is an eigenvector of $\rho (G)$ for a connected graph $G$, then $X$ is positive and unique.  The eigenequations for $A_{\alpha}(G)$ can be represented as the following form
\begin{eqnarray}
\rho(G)x_{u_i} = \alpha d_{u_i}x_{u_i} +(1 - \alpha) \sum_{u_iu_j \in E(G)} x_{u_j}.
\label{02}
\end{eqnarray}

Nikiforov et al. \cite{001} studied the $A_{\alpha}$-spectra of trees and determined the maximal  $A_{\alpha}$-spectral radius. It is known that a tree is a graph without cycles. If we replaced some vertices in a tree as a cycle, then this is an extension of the tree, that is, for an integer $k \geq 0,$
a cactus graph $C_n^k$ is a connected graph  such that  any two of its cycles have at most one common vertex.
Denoted by $\mathcal{C}_n^k$ be the set of all cacti  with $n$ vertices and $k$ cycles.  Let $C^c$ be a cactus graph in $\mathcal{C}_n^k$ such that all cycles (if any) have length $3$, that is $C^c$ contains $k$ cycles $vv_1v_1'v, vv_2v_2'v, \cdots,  vv_k v_k'v$ and $n-2k-1$ pendant edges $v u_1,vu_2, \cdots, vu_{n-2k-1}$. When $k = 0$, $C^c$ is a star. For other undefined notations and terminologies, refer to \cite{014}.
%\begin{figure}[htbp]
 %\centering
%\includegraphics[width=5in]{p1.png}
%\caption{ The cactus graph $C^c$.}
%\label{cc}
%\end{figure}

The cactus graph has been considered in mathematical literature, especially for the communication between graph theory and algebra. Borovi\'{c}anin and Petrovi\'{c} investigate  the properties of cacti with n vertices \cite{011}.
Chen and Zhou \cite{010} obtain the upper bound of the signless Laplacian spectral radius of cacti. Wu et al. \cite{016} find the spectral radius of cacti with $k$-pendant vertices. Shen et al. \cite{017} study the signless Laplacian spectral radius of cacti with given matching number. 

Inspired by the above results,
in this paper, we generalize the  $A_{\alpha}$-spectra from the trees to the cacti with $\alpha \in [0,1)$ and determine the largest $A_{\alpha}$-spectral radius  in $\mathcal{C}_n^k$. The extremal graph attaining  the sharp bound is proposed as well. Furthermore, we explore all eigenvalues of such extremal cacti. By using these outcomes, some previous results can be deduced, see\cite{001,010,011}.

\section{Main results and lemmas}
In this section, we first give some important lemmas that are used to our main proof.
The following lemma is another kind of results proposed by Nikiforov \cite{001} and Xue et al.\cite{002}. 
\begin{lem}\cite{001,002}
\label{lem1}
Let $A_{\alpha}(G)$ be the $A_{\alpha}$-matrix  of a connected graph $G$ with $0 \leq \alpha <1$,   $u \in S \subset V(G)$, and $v, w \in V (G)$ such that $S \subset N(v) \setminus (N(w) \cup \{w\})$.  Denote by $H$  the graph with vertex set $V(G)$ and edge set $E(G) \setminus \{uv, u \in S\} \cup \{uw, u \in S\}$,  and $X$  a unit eigenvector to
$\rho (A_{\alpha}(G))$. For $|S| \neq 0$, if  either \\$(i)$ $X^TA_{\alpha}(H)X \geq X^TA_{\alpha}(G)X$, 
or \\$(ii)$  $x_w \geq x_v$,
then $$\rho(H) > \rho(G).$$
\end{lem}

\begin{lem} 
\label{lem2}
Let $C_n^k$ be a cactus, $\alpha \in [0,1)$ and $C_l$ a cycle of $C_n^k$. If $\rho(C_n^k)$ is maximal, then $C_l$ is a triangle.
\end{lem}
\begin{proof} We prove it by a contradiction. Suppose that $C_n^k$ contains a cycle $C_l$ with the length $l \geq 4$. 

Let $uv$ be an edge in $C_l$ and $X$ be the unit eigenvector of $\rho(G)$.  Without loss of generality, assume that  $x_u \geq x_v$ and  $w \in V(C_l) \cap N(v) \setminus \{u\}$. We build a graph $H$ with vertex set $V(C_n^k)$ and edge set $E(C_n^k) \setminus \{vw\} \cup \{uw\}$. Then $H$ is a cactus graph and the length of $C_l$ decreases by $1$. By Lemma \ref{lem1}, we have $\rho(H) > \rho(C_n^k)$. This contradiction yields to our proof.
\end{proof}

\begin{lem} 
\label{lem3}\cite{002}
Let $G$ be a graph such that $u_0$ is a cut vertex, $u_0u_1 \in E(G)$ and a pendant path $u_1u_2 \cdots u_k$ is a component of $G- u_0$. For $\alpha \in [0,1)$,
if $X = (x_0, x_1,x_2, \cdots,x_k, \cdots, x_n)$ is a unit eigenvector of $\rho(G)$ corresponding to  the  vertex set $\{u_0,u_1,u_2, \cdots, u_k, \cdots, u_n\}$ and $\rho(G) > 2$, then $x_0>x_1>x_2>\cdots>x_k$.
\end{lem}

\begin{lem} 
\label{lem4}
Let $C_n^k$ be a cactus and $\alpha \in [0,1)$, if $\rho(C_n^k)$ is maximal,   the length of its pendant path   is $1$.
\end{lem}
\begin{proof}
We prove it by a contradiction. Suppose that there is a pendant path $u_1u_2\cdots u_k$ with $k \geq 2$ and $u_0$ is a cut vertex of degree at least $3$.  

Let $X = (x_0, x_1,x_2, \cdots, x_n)$ be a unit eigenvector of $G$ corresponding to $\rho(C_n^k) $  and vertex set $\{u_0,u_1,u_2, \cdots, u_n\}$.  Since $C_n^k$ is not a $2$-regular graph, then $\rho(C_n^k) > 2$. By Lemma \ref{lem3}, we have $x_0>x_1>x_2>\cdots>x_k$.

Let $H$ be a graph with vertex set $V(C_n^k)$ and edge set $E(C_n^k) \setminus \{u_1u_2\} \cup \{u_0u_2\}$. Then $H$ is a cactus graph.  Since $x_0 > x_1$, by Lemma \ref{lem1}, we have  $\rho(H) > \rho(C_n^k)$, which is a contradiction. We complete the proof.
\end{proof}

\begin{lem} 
\label{lem5}
Let $C_n^k$ be a cactus and $\alpha \in [0,1)$, if $\rho(C_n^k)$ is maximal, there is no proper cut edge.
\end{lem}
\begin{proof}
We prove it by a contradiction. Suppose that there exists a proper cut edge $uv$ such that $C_n^k - uv$ contains at least two components $G_1, G_2$ such that $|G_i| \geq 2$, $i=1, 2$. 

Let  $X$ be the unit eigenvector of $\rho(C_n^k)$. Without loss of generality, assume that $x_u \geq x_v$, $u \in V(G_1)$ and $v \in V(G_2)$. Let $S = N_{G}(v) \setminus \{u\}$. We set a new graph $H$ with vertex set $V(C_n^k)$ and edge set $E(C_n^k) \setminus \{vw, w\in S\} \cup \{uw, w \in S\}$. Then $H$ is a cactus graph. By Lemma \ref{lem1}, we have  $\rho(H) > \rho(C_n^k)$, which is a contradiction. The proof is completed.
\end{proof}

\subsection{The largest $A_{\alpha}$-spectral radius of cacti}
In this section, we provide the largest $A_{\alpha}$-spectral radius of a cactus graph $C_n^k$ in the set of cacti $\mathcal{C}_n^k$. 
\begin{thm}
Let $C_n^k  \in \mathcal{C}_n^k$ be a cactus and $\alpha \in [0,1)$. Then  $$\rho(C_n^k) \leq \rho(C^c).$$
\end{thm}
\begin{proof}
Let $\alpha \in [0, 1]$, and  $C_n^k$ be a cactus graph of order $n$ such that $\rho (A_{\alpha} (G))$ is maximal in $\mathcal{C}_n^k$. By Lemma \ref{lem2}, all cycles (if any) are of length 3.  By Lemma \ref{lem4}, all pendant paths are pendant edges. By Lemma \ref{lem5}, all  cycles are not connected by an edge or a path. 

Therefore, it suffices to prove that all cycles and pendant edges are sharing a common cut vertex. Next we prove the following claim.

\noindent {\bf Claim.} There exists a unique cut vertex in such $C_n^k$.

\noindent {\em Proof.} We prove it by a contradiction. Assume that there are at least two cut vertices $u, v$. By Lemma \ref{lem5}, $uv$ is not a cut edge. 

Let $N_u = \{w_u^1,w_u^2, \cdots, w_u^l\}$ and $N_v = \{w_v^1,w_v^2, \cdots, w_v^r\}$ be two neighborhoods of vertices $u$ and $v$. Without loss of generality, suppose that $x_u \geq x_v$ and $w_v^1$ has the shortest distance to the cut vertex $u$.   Denote $w_v^1, w_v^2$ and $v$  in a same cycle.  Now we build a new graph $H_1$ with vertex set $V(C_n^k)$ and edge set $E(C_n^k) \setminus \{v w_v^i, 3 \leq i \leq r\} \cup \{u w_v^i, 3 \leq i \leq r\}$. Note that the component number $w(H_1) = w(H) - 1$ and $H_1$ is still a cactus graph. By Lemma \ref{lem1}, we have  $\rho(H_1) > \rho(C_n^k)$. This is a contradiction that the chosen $C_n^k$ has the maximal $\rho$ in $\mathcal{C}_n^k$.

We can recursively apply the process using in Claim 1 and obtain the graph with the maximal $\rho$. Thus, we prove that the maximal $\rho$ attains  the cactus $C^c$.
\end{proof}

While we consider the relation between adjacent matrix $A(G)$, signless Laplacian matrix $Q(G)$, we can obtain the following corollary for the spectral radius $\rho_{A}$ and $\rho_{Q}$, respectively.
\begin{cor}\cite{011,010}
Let $C_n^k  \in \mathcal{C}_n^k$ be a cactus and $\alpha \in [0,1)$. Then  $$\rho(A(C_n^k)) \leq \rho(A(C^c)) \text{ and }\rho(Q(C_n^k)) \leq \rho(Q(C^c)).$$
\end{cor}

\subsection{The eigenvalues of $A_{\alpha}(C^c)$}
In this section, we determine the eigenvalues of $A_{\alpha}(C^c)$. Note that $2k+t+1 = n$. 
Let $I_n$ be the identity matrix of order $n.$
Let $J_n$  a matrix of all entries $1$ and $0_n$  a matrix of all entries $0$, respectively. 

\begin{thm} Label the vertices of $C^c$ as $v, v_1,v_2, \cdots, v_k, v_1',v_2' \cdots,v_k', u_1,u_2, \\ \cdots, u_t$ with $k, t \geq 0$. 
The eigenvalues of $A_{\alpha}(C^c)$ are $\alpha$,  $\alpha+1$(if $k \geq 2$, otherwise none),  $3\alpha - 1$ and the roots of $f(\lambda)=0$, where $ f(\lambda)=  (\alpha- \lambda)^3 + (n \alpha- 2\alpha+1) (\alpha- \lambda)^2 +[(1-n) \alpha^2 +(3n-4) \alpha +1-n](\alpha- \lambda) - t(1- \alpha)^2.$
\end{thm}

\begin{proof}
Since $A_{\alpha}(G) = \alpha D(G) +(1-\alpha)A(G)$, then the  $(i,j)$-entry of $ A_{\alpha}(C^c)$ is $d_{u_i} \alpha$ with $i=j$, $1- \alpha$ if $u_iu_j \in E(G)$ with $i \neq j$, and otherwise $0$.  Then
\begin{align} 
\setlength{\arraycolsep}{9pt}A_{\alpha} - \lambda I_n = 
\begin{bmatrix}
    (2k+t)\alpha -\lambda       & (1- \alpha)J_k^T & (1-\alpha)J_k^T &  (1- \alpha)J_t^T \\
    ~~\\
     (1- \alpha)J_k      & (2\alpha - \lambda)I_k   &   (1- \alpha)I_k     & 0 \\
     ~~\\
     (1- \alpha)J_k      & (1- \alpha)I_k       &   (2\alpha -\lambda) I_k     & 0\\~~\\
     (1- \alpha)J_t      & 0                           & 0                          &  (\alpha -\lambda) I_t
\end{bmatrix}.
%\nonumber
\end{align}

From the operations of this determinant $det [A_{\alpha} - \lambda I_n]$, we have

~~~~~~~~~~$det[A_{\alpha} - \lambda I_n]$
\[  =  
\left|\begin{array}{cccc} 
   (2k+t)\alpha -\lambda       & (1- \alpha)J_k^T & (1-\alpha)J_k^T &  (1- \alpha)J_t^T \\~~\\
     (1- \alpha)J_k      & (2\alpha - \lambda)I_k   &   (1- \alpha)I_k     & 0 \\~~\\
     (1- \alpha)J_k      & (1- \alpha)I_k       &   (2\alpha -\lambda) I_k     & 0\\~~\\
     (1- \alpha)J_t      & 0                           & 0                          &  (\alpha -\lambda) I_t
     \end{array}\right| 
\]
%~~~~~~~~~(~{\em Operations: Column $1$ $-$ (Column $i$) $\frac{1-\alpha}{\alpha-\lambda}$, $i \in [n-t+1, n]$}~)
\[  = (\alpha - \lambda )^t
\left|\begin{array}{cccc} 
    (2k+t)\alpha -\lambda- \frac{t(1-\alpha)^2}{\alpha - \lambda}       & (1- \alpha)J_k^T & (1-\alpha)J_k^T  \\~~\\
     (1- \alpha)J_k      & (2\alpha - \lambda)I_k   &   (1- \alpha)I_k      \\~~\\
     (1- \alpha)J_k      & (1- \alpha)I_k       &   (2\alpha -\lambda) I_k                   
     \end{array}\right| 
\]
%~~~~~~~~~(~{\em Operations: Column $j$ $-$ (Column $i$) $\frac{1-\alpha}{2\alpha-\lambda}$, $i \in [n-t-k+1, n-t], \\
%~~~~~~~~~~~~~~~~~~~~~~~~~~~j\in [1, n-t-k]$}~)
\[  =(\alpha - \lambda )^t
\left|\begin{array}{cccc} 
    (2k+t)\alpha -\lambda - \frac{t(1-\alpha)^2}{\alpha - \lambda}  & ((1- \alpha)- \frac{(1- \alpha)^2}{2\alpha - \lambda})J_k^T & (1- \alpha)J_k^T  \\
    - \frac{k(1- \alpha)^2}{2\alpha - \lambda}& \\~~\\
    ((1-\alpha)-\frac{(1- \alpha)^2}{2\alpha - \lambda})J_k  &  ((2\alpha - \lambda)- \frac{(1-\alpha)^2}{2\alpha - \lambda}) I_k   &   (1- \alpha)I_k      \\
    ~~\\
     0      & 0       &   (2\alpha -\lambda) I_k                   
     \end{array}\right| 
\]
%~~~~~~~~~(~{\em Operations: Column $1$ $-$ (Column $i$) $\frac{(1-\alpha)-\frac{(1- \alpha)^2}{2\alpha - \lambda}}{(2\alpha - \lambda)- \frac{(1-\alpha)^2}{2\alpha - \lambda}}$, $i \in [2, n-t-k]$}~)

\[  =(\alpha - \lambda )^t
\left|\begin{array}{cccc} 
    (2k+t)\alpha -\lambda - \frac{t(1-\alpha)^2}{\alpha - \lambda}- \frac{k(1- \alpha)^2}{2\alpha - \lambda} & ((1- \alpha)- \frac{(1- \alpha)^2}{2\alpha - \lambda})J_k^T & (1- \alpha)J_k^T  \\
   -  \frac{k(1-\alpha)^2(3\alpha-\lambda-1)}{(2\alpha-\lambda)(\alpha-\lambda+1)} &\\~~\\
  0  &  ((2\alpha - \lambda)- \frac{(1-\alpha)^2}{2\alpha - \lambda}) I_k   &   (1- \alpha)I_k      \\~~\\
     0      & 0       &   (2\alpha -\lambda) I_k                   
     \end{array}\right| 
\]

\[  =
(\alpha - \lambda )^t (2\alpha - \lambda)^k [\frac{(\alpha-\lambda+1)(3\alpha-\lambda-1)}{2\alpha-\lambda}]^k[(2k+t)\alpha -\lambda - \frac{t(1-\alpha)^2}{\alpha - \lambda}\]
\[ 
~~ - \frac{k(1- \alpha)^2}{2\alpha - \lambda} -
 \frac{k(1-\alpha)^2(3\alpha-\lambda-1)}{(2\alpha-\lambda)(\alpha-\lambda+1)}\]
\[  =
(\alpha - \lambda )^{t-1} (\alpha-\lambda+1)^{k-1}(3\alpha-\lambda-1)^k \{[(n-1)\alpha -\lambda](\alpha - \lambda ) (\alpha-\lambda+1) \]
\[ 
~~ - t(1-\alpha)^2(\alpha-\lambda+1)- 2k(1- \alpha)^2(\alpha - \lambda)\}. \]

In order to find the eigenvalues, we consider the characteristic equation  $$det[A_{\alpha} - \lambda I_n] = 0.$$ We have the roots   $\alpha$ of multiplicity $t-1$, $\alpha+1$ (if $k \geq 2$,  otherwise none) of multiplicity $k-1$, $3\alpha - 1$ of multiplicity $k$, and the other roots of $f(\lambda)=(n\alpha-\alpha -\lambda)(\alpha - \lambda ) (\alpha-\lambda+1)  - t(1-\alpha)^2(\alpha-\lambda+1)- 2k(1- \alpha)^2(\alpha - \lambda)=
 (\alpha- \lambda)^3 + (n \alpha- 2\alpha+1) (\alpha- \lambda)^2 +[(1-n) \alpha^2 +(3n-4) \alpha +1-n](\alpha- \lambda) - t(1- \alpha)^2
= 0.$  Therefore, these roots are the eigenvalues of $A_{\alpha}(C^c)$.
\end{proof}

Note that $n = 2k+t+1$ and the largest $A_{\alpha}$-spectral radius among trees attains at a star, that is $k= 0, t= n-1$. Applying such $k, t$ to $f(\lambda)$, we have the characteristic equation is
$$(\alpha - \lambda )^{n-2}[(n \alpha- \alpha -\lambda)(\alpha - \lambda) - (n-1)(1-\alpha)^2]= 0.$$ 
The  roots of this equation (or the eigenvalues of $A_{\alpha}$-matrix of  a star) are $\alpha$ of $n-2$ copies,  $\frac{\alpha n + \sqrt{\alpha^2n^2 + 4(n-1)(1-2\alpha)}}{2}$ and $\frac{\alpha n - \sqrt{\alpha^2n^2 + 4(n-1)(1-2\alpha)}}{2}$. Note that $\frac{\alpha n + \sqrt{\alpha^2n^2 + 4(n-1)(1-2\alpha)}}{2}$ is the largest one in these roots. In other words, we use a general method to prove the theorem  below. 

\begin{cor}\cite{001,005}
If $T$ is a tree with $n$ vertices and $0 \leq \alpha \leq 1$, then $$\rho(A_{\alpha}(T)) \leq \frac{\alpha n + \sqrt{\alpha^2n^2 + 4(n-1)(1-2\alpha)}}{2},$$ the equality holds if and only if $T$ is a star. In particular, the eigenvalues of $A_{\alpha}$-matrix of a star are $$\alpha, ~~\frac{\alpha n + \sqrt{\alpha^2n^2 + 4(n-1)(1-2\alpha)}}{2}  ~\text{and}~ \frac{\alpha n - \sqrt{\alpha^2n^2 + 4(n-1)(1-2\alpha)}}{2}.$$
\end{cor}

In addition, when $\alpha = 0$ or $\frac{1}{2}$,  the results of adjacent matrix from Lov\'{a}sz and  Pelik\'{a}n \cite{003}  and signless Laplacian matrix from Chen \cite{004} are deduced analogously, respectively.

{\bf Acknowledgement} The work was partially supported by the National Natural Science Foundation of China under Grants 11771172 and 11571134.

\end{document}